\documentclass{amsart}
\usepackage{amsmath}
\usepackage{amsthm}
\usepackage{amsfonts}
\usepackage{amscd}
\usepackage{mathrsfs}
\usepackage{amssymb}
\usepackage{euscript}
\usepackage{oldgerm}
\usepackage{graphicx}
\frakfamily

\theoremstyle{plain}
\newtheorem{thm}{Theorem}[section]
\newtheorem{prop}[thm]{Proposition}
\newtheorem{cor}[thm]{Corollary}
\newtheorem{lem}[thm]{Lemma}

\theoremstyle{definition}
\newtheorem{notation}[thm]{Notation}
\newtheorem{defn}[thm]{Definition}
\theoremstyle{remark}
\newtheorem{remark}[thm]{Remark}

\renewenvironment{proof}[1][Proof]{\textbf{#1.} }{\hfill \rule{0.5em}{0.5em}}
\numberwithin{equation}{section}

\title{Hypoelliptic heat kernel inequalities on Lie groups}
\author{Tai Melcher}
\address{Department of Mathematics, University
of Virginia, Charlottesville, VA, 22903}
\thanks{The author was supported in part by NSF Grants 99-71036 and DMS
0202939.} 
\email{melcher@virginia.edu}



\begin{document}

\begin{abstract}
This paper discusses the existence of gradient estimates for the heat kernel of a second order hypoelliptic operator on a manifold.  For elliptic operators, it is now standard that such estimates (satisfying certain conditions on coefficients) are equivalent to a lower bound on the Ricci tensor of the Riemannian metric.  For hypoelliptic operators, the associated ``Ricci curvature'' takes on the value $-\infty$ at points of degeneracy of the semi-Riemannian metric.  For this reason, the standard proofs for the elliptic theory fail in the hypoelliptic setting.

This paper presents recent results for hypoelliptic operators.  Malliavin calculus methods transfer the problem to one of determining certain infinite dimensional estimates.  Here, the underlying manifold is a Lie group, and the hypoelliptic operators are given by the sum of squares of left invariant vector fields.  In particular, ``$L^p$-type'' gradient estimates hold for $p\in(1,\infty)$, and the $p=2$ gradient estimate implies a Poincar\'e estimate in this context.
\end{abstract}

\maketitle
\setcounter{tocdepth}{1}
\tableofcontents

\section{Introduction}
\label{s:intro}

\subsection{Background}
Let $M$ be a manifold of dimension $d$, and let $\left\{  X_{i} \right\}_{i=1}^{k}$ be a set of smooth vector fields on $M$ satisfying
\begin{equation}
\label{e:HC}
T_{m}M=\mathrm{span}\left(  \{X(m):X\in\mathcal{L}\}\right) , \quad\forall~m\in M, \tag{HC}
\end{equation}
where $\mathcal{L}$ is the Lie algebra of vector fields generated by the collection $\left\{  X_{i}\right\}  _{i=1}^{k}$.  This assumption is the {\em H\"{o}rmander condition}, and the collection $\left\{  X_{i}\right\}  _{i=1}^{k}$ is a {\em H\"ormander set}.  Under this assumption, by a celebrated theorem of H\"{o}rmander, the operator 
\begin{equation}
\label{e:3}
L = \sum_{i=1}^{k}X_{i}^{2}
\end{equation}
is hypoelliptic.  Recall that a subelliptic operator $L$ is said to be {\em hypoelliptic} if $Lu\in C^\infty(\Omega)$ implies that $u\in C^\infty(\Omega)$, for all distributions $u\in C^\infty(\Omega)'$ on any open set $\Omega\subset_o M$.

\begin{notation}
\label{n:fxn}
Let $C^\infty_c(M)$ denote the set of smooth functions on $M$ with compact support, and let $C_b^\infty(M)$ denote the set of smooth, bounded functions on $M$.  When $M=\mathbb{R}^n$, let $C_p^\infty(\mathbb{R}^n)$ denote those functions $f\in C^\infty(\mathbb{R}^n)$ such that $f$ and all of its partial derivatives have at most polynomial growth.
\end{notation}

Let $\nabla=\left(  X_{1},\dots,X_{k}\right)$.  This paper continues the work begun in \cite{Melcher1}, considering $L^p$-type gradient inequalities of the form
\begin{equation}
\label{e:bl}
|\nabla e^{tL/2}f|^{p}\leq K_{p}(t)e^{tL/2}\left\vert \nabla f\right\vert^{p},\quad p\in\lbrack1,\infty),
\end{equation}
for $f\in C_{c}^{\infty}(M)$ and $t>0$.  For $p=1$, (\ref{e:bl}) is equivalent to a one parameter family of log Sobolev estimates for the heat kernel; for $p=2$, (\ref{e:bl}) is equivalent to a one parameter family of Poincar\'{e} estimates. The former has implications for hypercontractivity of an associated semigroup; see \cite{Gross75,Gross92}.  

When $L$ is an elliptic operator, a lower bound on the Ricci curvature is equivalent to the estimate (\ref{e:bl}) holding with some coefficients $K_p>0$ such that $K_p(0)=1$ and $\dot{K}_p(0)$ exists.  In particular, in the elliptic setting,  (\ref{e:bl}) holds with exponential coefficients $K_p(t)=e^{pkt}$, where $-2k$ is the lower bound on the Ricci curvature; see for example \cite{Bakry90,BakryEmery84,BakryEmery85}.  However, an operator $L$ of the form (\ref{e:3}) need not be elliptic. The principle symbol of $L$ at $\xi\in T_{m}^{\ast}M$ is given by $\sigma_{L}\left(  \xi\right)  =\sum_{i=1}^{k}\left[  \xi\left( X_{i}\right)  \right]  ^{2}$.  By definition, the operator $L$ is degenerate at points $m\in M$ where there exists $0\neq\xi\in T_{m}^{\ast}M$ such that $\sigma_{L}\left(  \xi\right)  =0$. At points of degeneracy of $L,$ the Ricci tensor is not well defined and should be interpreted to take the value $-\infty$ in some directions.  Thus, there exists no lower bound on the Ricci curvature in this case.  Nevertheless, it is reasonable to ask if inequalities of the form (\ref{e:bl}) might still hold, perhaps with some discontinuity in the coefficients $K_p$ near $t=0$.  In particular, under what conditions do functions $K_{p}\left(  t\right)  <\infty$ exist such that (\ref{e:bl}) is satisfied for all $f\in C_{c}^{\infty}(M)$ and $t>0$?  

The paper \cite{Melcher1} addressed the special case of the real three-dimensional Heisenberg Lie group, and the estimate (\ref{e:bl}) was proved to hold for all $p>1$ with a constant coefficient $K_p(t)\equiv K_p$, yielding a Poincar\'{e} estimate in this case.  Using analytic methods in \cite{Li06}, Li was able to prove (\ref{e:bl}) on the Heisenberg group for $p=1$, yielding the log Sobolev estimate.  Here in this paper, the case is addressed where the manifold $M$ is a general Lie group and the vector fields $\{X_i\}_{i=1}^k$ are invariant under left translation.

Related results appear in Kusuoka and Stroock \cite{KS3}, Picard \cite{Picard}, and Auscher, Coulhon, Duong, and Hofmann \cite{Coulhon}.  Also, \cite{Coulhon,Coulhon03} include some potential applications of the result proven here.

\subsection{Statement of results}

Let $G$ be a $d$-dimensional Lie group with Lie algebra $\mathfrak{g}=\mathrm{Lie}(G)$ and identity element $e$.  Let $L_g$ denote left translation by an element $g\in G$, and let $R_g$ denote right 
translation.  Suppose $\{X_{i}\}_{i=1}^{k}\subset\mathfrak{g}$ is a linearly independent Lie generating set; that is, there exists some $m\in\mathbb{N}$ such that
\begin{multline}
\label{e:gen} 
\mathrm{span}\big\{X_i,[X_{i_1},X_{i_2}],[X_{i_1},[X_{i_2},X_{i_3}]],\ldots,
   [X_{i_1},[\cdots,[X_{i_{m-1}},X_{i_m}]\cdots]] :  \\ 
   i,i_r \in\{1,\ldots,k\}, r\in\{1,\ldots,m\} \big\} = \mathfrak{g}.
\end{multline}

\begin{notation}
\label{n:Sigma}
Let $\Sigma=\Sigma_0:=\{X_1,\ldots,X_k\}$ and $\Sigma_r$ be defined inductively by
\[ \Sigma_r:=\{[X_i,V]:V\in\Sigma_{r-1},i=1,\ldots,k\}, \]
for all $r\in\mathbb{N}$.  Since $\{X_i\}_{i=1}^k$ is a Lie generating set, there is a finite $m$
such that 
\[   \mathrm{span}\left(\cup_{r=0}^m \Sigma_r \right) = \mathfrak{g}. \]
Let $\mathfrak{g}_{0}:=\mathrm{span}(\Sigma_0)$, and let $\{Y_j\}_{j=1}^{d-k}\subset\cup_{r=1}^m \Sigma_r$ be a basis of $\mathfrak{g}/\mathfrak{g}_0$.  Define an inner product $\left\langle  \cdot,\cdot \right\rangle$ on $\mathfrak{g}$ by making $\{X_i\}_{i=1}^k\cup\{Y_j\}_{j=1}^{d-k}$ an orthonormal set.  Note then that $\{X_i\}_{i=1}^k$ is an orthonormal basis of $\mathfrak{g}_0$.  Extend $\left\langle  \cdot
,\cdot\right\rangle  $ to a right invariant metric on $G$ by defining
$\left\langle \cdot,\cdot \right\rangle_g : T_gG \times T_gG \rightarrow \mathbb{R}$ as
\[
\left\langle v,w \right\rangle  _g := \left\langle  R_{g^{-1}*}v , R_{g^{-1}*}w\right\rangle  ,\quad\text{ for all }v,w\in T_gG.
\]
The $g$ subscript will be suppressed when there is no chance of confusion. 
\end{notation}

\begin{notation}
Given an element $X\in\mathfrak{g}$, let $\tilde{X}$ denote the left invariant
vector field on $G$ such that $\tilde{X}(e)=X$, where $e$ is the identity of $G$.
Recall that $\tilde{X}$ left invariant means that the vector field commutes
with left translation in the following way:
\[ \tilde{X}(f\circ L_g) = (\tilde{X} f) \circ L_g, \]
for all $f\in C^1(G)$.
Similarly, let $\hat{X}$ denote the right invariant vector field associated to $X$.
\end{notation}


\begin{defn}
The left invariant {\em gradient} on $G$ is the operator on $C^1(G)$ given by
\[ \nabla := (\tilde{X}_1,\ldots,\tilde{X}_k). \]
The {\em subLaplacian} on $G$ is the second-order operator acting on $C^2(G)$ given by
\[  L:= \sum_{i=1}^k \tilde{X}_i^2. \]
\end{defn}

\begin{remark}
\label{r:Lhypo}
Since $\{X_i\}_{i=1}^k$ is a Lie generating set, $\{\tilde{X}_i\}_{i=1}^k$ satisfies the H\"{o}rmander condition (\ref{e:HC}) and H\"{o}rmander's theorem \cite{Hormander67} implies that $L$ is hypoelliptic.
\end{remark}

Let $L^2(G)$ denote the space of square integrable functions on $G$ with respect to right invariant Haar measure. 
Then $L$ is a densely defined, symmetric operator on $L^2(G)$ and the symmetric bilinear form associated to $L$ is given by $\mathcal{E}^0(f_1,f_2) := (-Lf_1,f_2)_{L^2(G)}$.  Note that $\mathcal{E}^0$ is positive, and so
$\mathcal{E}^0$ is closable.  The minimal closure $\mathcal{E}$ is associated to a self-adjoint operator $\bar{L}$ which is an extension of $L$, called the Friedrichs extension of $L$.  

\begin{defn}
\label{d:Pt}
Let $P_t$ denote the {\em heat semigroup} $e^{t\bar{L}/2}$, where $\bar{L}$ is the Friedrichs extension of $L|_{C_c^\infty(G)}$ to $L^2(G,dg)$ with $dg$ right Haar measure on $G$.  By the left invariance of $L$ and the satisfaction of the H\"ormander condition, $P_t$ admits a left convolution kernel $p_t$ such that 
\begin{equation*}
P_tf(h) =f *p_t(h) = \int_G f(hg)p_t(g) \,dg, 
\end{equation*}
for all $f\in C_c^\infty(G)$.  The function $p_t$ is called the {\em heat kernel} of $G$. 
\end{defn}

The operator $P_t$ is a symmetric Markov semigroup.  By Remark \ref{r:Lhypo}, $L$ is a hypoelliptic operator, and so $p_t$ is a smooth density on $G$.   In the sequel, let $L$ denote its own Friedrichs extension.  For the standard semigroup theory used here, see for example \cite{Davies}. 

\begin{notation}
\label{n:best}
Let $K_p(t)  $ be the best function such that
\begin{equation}
\label{e:MAIN}
|\nabla P_{t}f|^{p}\le K_{p}(t)  P_{t}|\nabla f|^{p}, \quad  p\in[1,\infty), \tag{$I_p$}
\end{equation}
for all $f\in C_{c}^{\infty}(G)$ and $t>0$.
\end{notation}

\begin{thm}
\label{t:main}
For all $p\in(1,\infty)$, $K_p(t)<\infty$ for all $t>0$.  If $G$ is a nilpotent Lie group, then there exists a constant $K_p<\infty$ such that $K_p(t)\le K_p$ for all $t>0$.  
\end{thm}

This theorem was established in \cite{Melcher1} in the case of the real three-dimensional Heisenberg Lie group.  The method of proof in this case is analogous.  The heat kernel $p_{t}(g)\,dg$ may be realized as the distribution in $t$ of the Cartan rolling map on $G$, the process $\xi$ satisfying the Stratonovich stochastic differential equation
\[ d\xi_t = \sum_{i=1}^k \tilde{X}_i(\xi_t) \circ db^i_t, \text{ with } \xi_0 = e, \]
where $b^1,\ldots,b^k$ are $k$ independent real-valued Brownian motions.  Thus, for all $f\in C_c^\infty(G)$, 
\[ P_tf(e)=\mathbb{E}[f(\xi_t)]. \]  
Sections \ref{s:SDE} and \ref{s:covariance} discuss properties of $\xi$.  This representation of $P_{t}$ transforms the finite dimensional problem to a problem on Wiener space.    Section \ref{s:lift} describes a standard ``lifting" procedure which constructs vector fields $\mathbf{X}_i$ on Wiener space from the vector fields $\tilde{X}_i$ via the map $\xi$. Then Malliavin's probabilistic techniques on proving hypoellipticity give componentwise bounds of $P_t(\tilde{X}_i f)(e) = \mathbb{E}[(\tilde{X}_if)(\xi_t)]=\mathbb{E}[\mathbf{X}_i(f(\xi_t))]$.  Section \ref{s:Mall.def} reviews some calculus on Wiener space necessary for this argument.

Section \ref{s:Lie} contains the proof of Theorem {\ref{t:main}.  Results from Section \ref{s:Malliavin} show that for a Lie group $G$, $K_p(t)<\infty$ for all $t>0$;  however, this method does not give any estimates on the behavior of $K_p$ with respect to $t$.  In a generalization of the Heisenberg scaling argument in \cite{Melcher1}, Section \ref{s:poincare} addresses the special case of nilpotent and stratified groups.  When $G$ is stratified, dilation arguments imply that the coefficients $K_{p}$ are independent of the $t$ parameter.  When $G$ is nilpotent, covering $G$ with a stratified group shows that there is a constant $K_p$ such that $K_p(t)\le K_p$ for all $t>0$, and this completes the proof of Theorem \ref{t:main}.  This implies the following Poincar\'e estimate for the heat kernel measure in this context. 

\begin{thm}
\label{t.poinc1} 
Suppose $G$ is a nilpotent Lie group with identity element $e$.  Then
\[
P_{t}f^{2}(e)-(P_{t}f)^{2}(e)\leq K_{2}tP_{t}|\nabla f|^{2}(e),
\]
for all $f\in C_c^{\infty}(G)$ and $t>0$, where $K_2$ is the constant in Theorem \ref{t:main} for $p=2$. 
\end{thm}

Note that this theorem gives an improvement in the elliptic case with negative curvature, giving linear coefficients where the estimate was previously known only with coefficients of exponential growth.  This is stated explicitly in Corollary \ref{c:elliptic}.  It could be conjectured that this is true for every Riemannian manifold; that is, for any Riemannian manifold equipped with a Laplace Beltrami operator, Poincar\'e estimates for the associated heat kernel hold with linear coefficients.

\subsection*{Acknowledgement.} I thank Bruce Driver for suggesting
this problem and for many valuable discussions throughout the
preparation of this work.

\section{Wiener calculus over $G$}
\label{s:Malliavin}

\subsection{Review of calculus on Wiener space}
\label{s:Mall.def}

This section contains a brief introduction to basic Wiener space definitions and notions of differentiability.  For a more complete exposition, consult \cite{Driver04,IkedaWatanabe89,Nualart95} and references contained therein.

Let $(\mathscr{W}(\mathbb{R}^{k}),\mathcal{F},\mu)$ denote classical k-dimensional
Wiener space. That is, $\mathscr{W}=\mathscr{W}(\mathbb{R}^k)$ is the Banach space of continuous paths $\omega:[0,1]\rightarrow\mathbb{R}^{k}$ such that $\omega_0=0$, equipped with
the supremum norm
\[
\Vert\omega\Vert=\max_{t\in\lbrack0,1]}|\omega_t|,
\]
$\mu$ is standard Wiener measure, and $\mathcal{F}$ is the completion of the
Borel $\sigma$-field on $\mathscr{W}$ with respect to $\mu$.  By definition of $\mu,$ the process
\[
b_t(\omega) = (b_t^1(\omega), \ldots , b_t^k(\omega) ) = \omega_t
\]
is an $\mathbb{R}^k$ Brownian motion. For those $\omega\in \mathscr{W}$ which are
absolutely continuous, let
\[
E(\omega):=\int_{0}^{1}|\dot{\omega}_s|^{2}\,ds
\]
denote the energy of $\omega$. The Cameron-Martin space is the Hilbert space
of finite energy paths,
\[
\mathscr{H}=\mathscr{H}(\mathbb{R}^k):=\{\omega\in \mathscr{W}(\mathbb{R}^k):\omega\text{ is
absolutely continuous and }E(\omega)<\infty\},
\]
equipped with the inner product
\[
(h,k)_{\mathscr{H}}:=\int_{0}^{1} \dot{h}_s\cdot\dot{k}_s\,ds, \quad\text{ for all }h,k\in
\mathscr{H}.
\]
More generally, for any finite dimensional vector space $V$ equipped with an inner product, let $\mathscr{W}(V)$  denote path space on $V$, and $\mathscr{H}(V)$ denote the set of Cameron-Martin paths, where the definitions are completely analogous, replacing the inner products and norms where necessary.

\begin{defn}
Denote by $\mathcal{S}$ the class of {\em smooth cylinder functionals},
random variables $F:\mathscr{W}\rightarrow\mathbb{R}$ such that
\begin{equation}
\label{e:cyl}
F(\omega)=f(\omega_{t_1},\ldots,\omega_{t_n}),
\end{equation}
for some $n\geq1,$ $0 < t_{1} < \cdots < t_{n} \le 1$, and function $f\in
C^\infty_{p}((\mathbb{R}^{k})^n)$ (see Notation \ref{n:fxn}).  For $E$ be a real separable Hilbert space, let $\mathcal{S}_{E}$ be the set of $E$-valued smooth cylinder functions $F:\mathscr{W}\rightarrow E$ of the form
\begin{equation}
\label{e:Ecyl}
F=\sum_{j=1}^{m} F_{j} e_{j},
\end{equation}
for some $m\ge1$, $e_{j}\in E$, and $F_{j}\in\mathcal{S}$.
\end{defn}

\begin{defn}
Fix $h\in\mathscr{H}$.  The directional derivative of a smooth cylinder functional $F\in\mathcal{S}$ of the form (\ref{e:cyl}) along $h$ is given by
\begin{equation*}
\partial_{h}F(\omega):=\frac{d}{d\epsilon}\bigg|_{0}F(\omega+\epsilon h) = 
  \sum_{i=1}^n  \nabla^i f(\omega_{t_1},\ldots,\omega_{t_n})\cdot h_{t_i},
\end{equation*}
where $\nabla^i f$ is the gradient of $f$ with respect to the $i^{th}$ variable.
\end{defn}

The following integration by parts result is standard; see for example Theorem 8.2.2 of Hsu \cite{Hsu03}.
\begin{prop}
\label{p:dh*}
Let $F,G\in\mathcal{S}$ and $h\in\mathscr{H}$.  Then
\[ (\partial_h F,G)_\mathscr{H} = (F,\partial_h^* G)_\mathscr{H}, \]
where $\partial_h^*=-\partial_h + \int_0^1 \dot{h}_s\cdot db_s$.
\end{prop}

\begin{defn}
The gradient of a smooth cylinder functional $F\in\mathcal{S}$ is the random process $D_tF$
taking values in $\mathscr{H}$ such that $(DF,h)_{\mathscr{H}}=\partial_{h}F$.  It may be determined that, for $F$ of the form (\ref{e:cyl}),
\[
D_{t}F = \sum_{i=1}^{n} \nabla^if(\omega_{t_1},\ldots,\omega_{t_n}) (t_i\wedge t),
\]
where $s\wedge t=\min\{s,t\}$.  For $F\in\mathcal{S}_E$ of the form (\ref{e:Ecyl}), define the derivative $D_tF$ to be the random process taking values in $\mathscr{H}\otimes E$ given by
\[ D_tF := \sum_{j=1}^m D_tF_j\otimes e_j. \]
\end{defn}

Iterations of the derivative for smooth functionals $F\in\mathcal{S}$ are given by
\[
D_{t_{1},\ldots,t_{k}}^{k}F=D_{t_{1}}\cdots D_{t_{k}}F \in \mathscr{H}^{\otimes k},
\]
for $k\in\mathbb{N}$.  For $F\in\mathcal{S}_E$,
\[
D^{k} F = \sum_{j=1}^{m} D^{k} F_{j} \otimes e_{j},
\]
and these are measurable functions defined almost everywhere on $[0,1]^{k}\times \mathscr{W}$.  The operator $D$ on $\mathcal{S}_E$ is closable, and there exist closed extensions $D^{k}$ to $L^{p}(\mathscr{W},\mathscr{H}^{\otimes k}\otimes E)$; see, for example \cite{Nualart95}, Theorem 8.28 of \cite{Hsu03}, or Theorem 8.5 of \cite{IkedaWatanabe89}.
Denote the closure of the derivative operator also by $D$ and the domain of $D^{k}$ in $L^{p}([0,1]^{k}\times \mathscr{W})$ by ${\mathcal{D}}^{k,p}$, which is the completion of the family of smooth Wiener functionals $\mathcal{S}$ with respect to the seminorm $\Vert\cdot\Vert_{k,p,E}$ on $\mathcal{S}_E$ given by
\[
\|F\|_{k,p,E} := \left( \sum_{j=0}^{k} \mathbb{E}(\|D^{j}F\|
   _{\mathscr{H}^{\otimes j}\otimes E}^{p}) \right)^{1/p},
\]
for any $p\ge1$.  Let
\[
\mathcal{D}^{k,\infty}(E) := \bigcap_{p>1} \mathcal{D}^{k,p}(E) \text{ and } 
   \mathcal{D}^{\infty}(E) := \bigcap_{p>1}\bigcap_{k\ge1}{\mathcal{D}}^{k,p}(E).
\]
When $E=\mathbb{R}$, write $\mathcal{D}^{k,p}(\mathbb{R})=\mathcal{D}^{k,p}$, $\mathcal{D}^{k,\infty}(\mathbb{R})=\mathcal{D}^{k,\infty}$, and $\mathcal{D}^\infty(\mathbb{R})=\mathcal{D}^\infty$.

\begin{defn}
Let $D^{*}$ denote the $L^2(\mu)$-adjoint of the derivative operator $D$, which has
domain in $L^{2}(\mathscr{W}\times[0,1],\mathscr{H})$ consisting of functions $G$ such that
\[
| \mathbb{E}[(DF,G)_{\mathscr{H}}]| \le C\|F\|_{L^{2}(\mu)},
\]
for all $F\in{\mathcal{D}}^{1,2}$, where $C$ is a constant depending on $G$.
For those functions $G$ in the domain of $D^{*}$, $D^{*}G$ is the element of
$L^{2}(\mu)$ such that
\[
\mathbb{E}[FD^{*}G] = \mathbb{E}[(DF,G)_{\mathscr{H}}].
\]
\end{defn}

It is known that $D$ is a continuous operator from $\mathcal{D}^{\infty}$ to
$\mathcal{D}^{\infty}(\mathscr{H})$, and similarly, $D^{\ast}$ is continuous from
$\mathcal{D}^{\infty}(\mathscr{H})$ to $\mathcal{D}^{\infty}$; see for example
Theorem V-8.1 and its corollary in \cite{IkedaWatanabe89} . 

Malliavin \cite{Malliavin78b,Malliavin78c} introduced the notion of derivatives of Wiener functionals and applied it to the regularity of probability laws induced by the solutions to stochastic differential equations at fixed times.  The notion of Sobolev spaces of Wiener functionals was first introduced by Shigekawa \cite{Shigekawa80} and Stroock \cite{Stroock80a,Stroock80b}.

\subsection{Rolling map}
\label{s:SDE}

Now, let $G$ be a Lie group with identity $e$ and Lie algebra $\mathrm{Lie}(G)=\mathfrak{g}$, and suppose $\{X_{i}\}_{i=1}^{k}\subset\mathfrak{g}$ is a linearly independent Lie generating set, in the sense of Equation (\ref{e:gen}).  Recall that $\{X_i\}_{i=1}^k$ is an orthonormal basis of the subspace $\mathfrak{g}_0=\mathrm{span}(\{X_i\}_{i=1}^k)$ with respect to the inner product defined on $\mathfrak{g}$.

\begin{notation}
Let $\operatorname{Ad}:G\rightarrow\mathrm{End}(\mathfrak{g})$ denote the adjoint representation of $G$ with differential $\operatorname{ad}:= d(\operatorname{Ad}):\mathfrak{g}\rightarrow\mathrm{End}(\mathfrak{g})$.  That is, $\operatorname{Ad}(g)=\operatorname{Ad}_g=L_{g*}R_{g^{-1}*}$, for all $g\in G$, and $\operatorname{ad}(X)=\operatorname{ad}_X=[X,\cdot]$, for all $X\in\mathfrak{g}$.  For any function $\varphi\in C^{1}(G)$, define $\hat{\nabla}\varphi,\tilde{\nabla}\varphi:G\rightarrow\mathfrak{g}$ such that, for any $g\in G$ and $X\in\mathfrak{g}$,
\begin{equation*}
\begin{split}
 \left\langle \hat{\nabla}\varphi(g),X \right\rangle &:= \left\langle d\varphi(g),R_{g*}X \right\rangle = (\hat{X}\varphi)(g) \\
\left\langle \tilde{\nabla}\varphi(g),X\right\rangle &:= \left\langle d\varphi(g),L_{g*}X\right\rangle  =(\tilde{X}\varphi)(g). 
\end{split}
\end{equation*}
\end{notation}

The sequel will use the following facts:
\begin{equation}
\label{e:d2.0}
\begin{split}
\left\langle  \hat{\nabla}\varphi(g),X\right\rangle  
   &= \left\langle d\varphi(g),L_{g*}L_{g^{-1}*}R_{g*}X \right\rangle  \\
   &=\left\langle d\varphi(g),L_{g*}\operatorname{Ad}_{g^{-1}}X \right\rangle  
   =\left\langle  \tilde{\nabla}\varphi(g) , \operatorname{Ad}_{g^{-1}}X \right\rangle 
\end{split}
\end{equation}
and similarly
\begin{equation}
\label{e:d2.1}
\left\langle  \tilde{\nabla}\varphi(g) , X \right\rangle = \left\langle \hat{\nabla}\varphi(g) , \operatorname{Ad}_gX \right\rangle  .
\end{equation}

Now suppose $\{b_t^i\}_{i=1}^k$ are $k$ independent real-valued Brownian motions.  Then
\[ \vec{b}_t := X_i b_t^i := \sum_{i=1}^k X_i b^i_t \]
is a $(\mathfrak{g}_0,\langle\cdot,\cdot\rangle)$ Brownian motion.  In the sequel, the convention of summing over repeated upper and lower indices will be observed.  Let $\xi:[0,1]\times \mathscr{W}\rightarrow G$ denote the solution to the Stratonovich stochastic differential equation
\begin{equation}
\label{e:SDE}
d\xi_{t} = \xi_t \circ d\vec{b}_t := L_{\xi_t*}\circ d\vec{b}_t 
   = L_{\xi_t*}X_i \circ db_t^i = \tilde{X}_i(\xi_t) \circ db_t^i,
\text{ with } \xi_0 = e.
\end{equation}
The solution $\xi$ exists by the standard theory;  see, for example, Theorem V-1.1 of \cite{IkedaWatanabe89}.  Additionally, Remark V-10.3 of \cite{IkedaWatanabe89} implies that $P_{t}=e^{tL/2}$, with $L=\sum_{i=1}^k \tilde{X}_i^2$, is the associated Markov diffusion semigroup to $\xi$, where $P_t$ is as defined in Definition \ref{d:Pt}; that is, $\nu_{t}:=(\xi_{t})_*\mu=p_{t}(g)\,dg$ is the density of the transition probability of the diffusion process $\xi_{t}$, where $dg$ denotes right Haar measure, and
\begin{equation}
\label{e:integral}
(P_{t}f)(e)=\mathbb{E}[f(\xi_{t})],
\end{equation}
for any $f\in C_c^{\infty}(G)$, where the right hand side is expectation
conditioned on $\xi_{0}=e.$

The following theorem is proved in \cite{Melcher3}.
\begin{thm}
\label{t:xi.smooth}
For any $f\in C_c^\infty(G)$, $f(\xi_t)\in\mathcal{D}^\infty$ for all $t\in[0,1]$.  In particular, $D[f(\xi_t)]\in\mathscr{H}\otimes\mathbb{R}^k$ and 
\begin{equation}
\label{e:Dxi}
(D[f(\xi_t)])^i = \left\langle \hat{\nabla}f(\xi_t),  \int_0^{\cdot\wedge t}\operatorname{Ad}_{\xi_\tau} X_i \,d\tau\right\rangle,
\end{equation}
for $i=1,\ldots,k$, componentwise in $\mathscr{H}$, and, for any $h\in\mathscr{H}$,
\begin{equation} 
\label{e:dhf} 
\begin{split}
\partial_h f(\xi_t) 
	= \left\langle \hat{\nabla}f(\xi_t),\int_0^t \operatorname{Ad}_{\xi_s} X_i\dot{h}_s^i\,ds \right\rangle
	= \left\langle df(\xi_t), R_{\xi_t*}\int_0^t \operatorname{Ad}_{\xi_s} X_i\dot{h}_s^i\,ds \right\rangle.
\end{split}
\end{equation}
\end{thm}

\begin{notation}
Let $\bigcap_{p>1}L^{p}(\mu)=:L^{\infty-}(\mu)$.
\end{notation}

\subsection{Covariance matrix}
\label{s:covariance}

The Malliavin covariance matrix of $\xi$ is the matrix
$\sigma_t(\omega):=\xi_t'(\omega)\xi_t'(\omega)^*:T_{\xi_t(\omega)}G\rightarrow T_{\xi_t(\omega)}G$,
where $\xi_t'(\omega):\mathscr{H}\rightarrow T_{\xi_t(\omega)}G$ is the Frechet derivative given by
\[ \xi_t'(\omega)h := \frac{d}{d\epsilon}\bigg|_0 \xi_t(\omega+\epsilon h), \] 
for all $h\in \mathscr{H}$, and its adjoint $\xi_t'(\omega)^*:T_{\xi_t(\omega)}\rightarrow \mathscr{H}$ is computed relative to the Cameron-Martin inner product on $\mathscr{H}$ and the chosen metric on $G$.   
Note that Equation (\ref{e:dhf}) implies that
\begin{equation}
\label{e:dxih} 
\xi_t'(\omega) h = R_{\xi_t*}\int_0^t \operatorname{Ad}_{\xi_s} X_i\dot{h}^i_s\,ds
\end{equation}

\begin{notation}
In the following, let $\operatorname{Ad}_{\xi_t}^{^{\dag}}$ denote the adjoint of $\operatorname{Ad}_{\xi_t}$ as an operator on $\mathfrak{g}$, and let $P:\mathfrak{g}\rightarrow\mathfrak{g}_0$ be orthogonal projection onto the subspace $\mathfrak{g}_0$. 
\end{notation}

\begin{thm}
\label{t:inv.smooth}
The Malliavin covariance matrix of $\xi$ is
\begin{equation}
\label{e:sigma} 
\sigma_t
   := \xi_t'(\omega)\xi_t'(\omega)^*
   = R_{\xi_t*}\left(\int_0^t \operatorname{Ad}_{\xi_s}P\operatorname{Ad}_{\xi_s}^\dag\,ds \right)R_{\xi_t*}^{\operatorname{tr}}.
\end{equation}
Let $\bar{\sigma}_t = \int_0^t \operatorname{Ad}_{\xi_s}P\,\operatorname{Ad}_{\xi_s}^\dag\,ds$, and $\Delta_t:=\det\bar{\sigma}_t$.  Then $\Delta_t>0$ a.e., and so $\bar{\sigma}_{t}$
is invertible a.e. for $t>0$.  Moreover, 
\[
\Delta_t^{-1}\in L^{\infty-}(\mu).
\]
\end{thm}

\begin{proof}
To determine $\sigma_t= \xi_t'(\omega)\xi_t'(\omega)^*$, first compute $\xi_t'(\omega)^*:T_{\xi_t(\omega)}G \rightarrow \mathscr{H}$, the adjoint in $\xi_t'(\omega)$ with respect to the Cameron-Martin inner product and the right invariant metric on $TG$.  By Equation (\ref{e:dxih}), for any $X\in\mathfrak{g}$,
\begin{align*}
\left(  \xi_{t}^{\prime}(\omega)^{\ast}( R_{\xi_t*}X)  ,h\right)_{\mathscr{H}}    
   &=\left\langle  R_{\xi_t*}X  ,\xi_{t}^{\prime}(\omega)h\right\rangle \\
   &=\left\langle R_{\xi_t*}X  , R_{\xi_t*}  \int_{0}^{t}\operatorname{Ad}_{\xi_s}X_i\dot{h}^i_{s} \,ds \right\rangle \\
   &=\left\langle X,\int_{0}^{t}\operatorname{Ad}_{\xi_{s}}X_i\dot{h}^i_{s} \, ds\right\rangle 
    =\int_{0}^{t}\left\langle \operatorname{Ad}_{\xi_{s}}^{^{\dag}}X,X_i\right\rangle\dot{h}^i_{s} ds,
\end{align*}
where the penultimate equality follows from the right invariance of the 
metric on $G$.  It then follows that
\begin{equation}
\label{e:xi'*}
\frac{d}{ds} \left[  \xi_{t}^{\prime}(\omega)^{\ast}(R_{\xi_t*}X)  \right]^i_{s}
   = 1_{s\le t} \left\langle \operatorname{Ad}_{\xi_s}^\dag X, X_i \right\rangle,
\end{equation}
componentwise in $\mathscr{H}$.  Combining Equations \eqref{e:dxih} and \eqref{e:xi'*},
\begin{align*} 
\xi_t'(\omega)\xi_t'(\omega)^*(R_{\xi_t*}X) 
   &= R_{\xi_t*}\int_0^t \operatorname{Ad}_{\xi_s} X_i\frac{d}{ds} \left[ 
   	\xi_{t}^{\prime}(\omega)^{\ast}(R_{\xi_t*}X)  \right]^i_{s}\,ds \\
   &= \sum_{i=1}^k R_{\xi_t*} \int_0^t \operatorname{Ad}_{\xi_s} X_i \left\langle\operatorname{Ad}_{\xi_s}^\dag X,X_i \right\rangle \,ds \\
   &= R_{\xi_t*} \int_0^t \operatorname{Ad}_{\xi_s}P\operatorname{Ad}_{\xi_s}^\dag X ds, 
\end{align*}
and Equation (\ref{e:sigma}) follows.

The proof that $\Delta_t>0$ and $\Delta_t^{-1}\in L^{\infty-}(\mu)$ is by now standard and relies on satisfaction of the H\"{o}rmander bracket condition, $\mathrm{Lie}(\{X_i\}_{i=1}^k) = \mathfrak{g}$;  for example, a simple adaptation of the proof of Theorem 8.6 in Driver \cite{Driver04} will work.
\end{proof}

\begin{remark}
\label{rem:mall}  
By the general theory, Theorem \ref{t:inv.smooth} implies $\nu_{t}=\mathrm{Law}(\xi_{t})$ is a smooth measure; see for example Remark V-10.3 of \cite{IkedaWatanabe89}.
\end{remark}

\subsection{Lifted vector fields and $L^2$-adjoints}
\label{s:lift}

Throughout this section, $t\in[0,1]$ will be fixed.
\begin{defn}
\label{d:lift}
Given $X\in \mathfrak{g}$, let $\tilde{X}$ be the associated left invariant vector field on $G$.
Define the ``{\em lifted vector field}" $\mathbf{X}$ of $\tilde{X}$ as
\begin{equation}
\label{e:lift}
\mathbf{X} = \mathbf{X}^t
   :=  \xi_{t}^{\prime}(\omega)^{\ast}\left[\xi_t'(\omega)\xi_t'(\omega)^*\right] ^{-1}\tilde{X}(\xi_{t}) 
   = \xi_{t}^{\prime}(\omega)^{\ast}\sigma_t^{-1}\tilde{X}(\xi_{t}) \in \mathscr{H}, 
\end{equation}
acting on functions $F\in\mathcal{D}^{1,2}$ by
\[ \mathbf{X} F = (DF,\mathbf{X})_\mathscr{H}. \]
\end{defn}

\begin{prop}
\label{p:A.smooth}  
For any $X\in\mathfrak{g}$, $\mathbf{X}\in\mathcal{D}^\infty(\mathscr{H})$, and  
\[ \mathbf{X}[f(\xi_t)] = (\tilde{X} f)(\xi_t), \]
for any $f\in C^\infty(G)$, 
\end{prop}

\begin{proof}
Combining Equations (\ref{e:sigma}) and (\ref{e:xi'*}) gives
\[
\frac{d}{ds} \mathbf{X}^i_s = 1_{s\le t} \left\langle \operatorname{Ad}_{\xi_s}^\dagger \left( \int_0^t \operatorname{Ad}_{\xi_r}
     P\,\operatorname{Ad}_{\xi_r}^\dagger \,dr \right)^{-1} \operatorname{Ad}_{\xi_t}X, X_i\right\rangle.
\]
Thus, rewrite Equation (\ref{e:lift}) explicitly as
\begin{align}
\mathbf{X}^i &= \int_0^{\cdot\wedge t} \left\langle \operatorname{Ad}_{\xi_s}^\dagger \left( \int_0^t \operatorname{Ad}_{\xi_r}
     P\,\operatorname{Ad}_{\xi_r}^\dagger \,dr \right)^{-1} \operatorname{Ad}_{\xi_t}X, X_i \right\rangle \,ds \notag \\
\label{e:Xs}
   &= \left\langle \left(\int_0^{\cdot\wedge t} \operatorname{Ad}_{\xi_s}^\dag\, ds \right) \left( \int_0^t \operatorname{Ad}_{\xi_s}
     P\,\operatorname{Ad}_{\xi_s}^\dagger \,ds \right)^{-1} \operatorname{Ad}_{\xi_t}X, X_i \right\rangle.
\end{align}
A standard argument shows that $\overline{W}=\int_0^\cdot \operatorname{Ad}_{\xi_s} \,ds \in\mathcal{D}^\infty(\mathscr{H}(\mathrm{End}(\mathfrak{g})))$; see for example Proposition 5 of \cite{app}.  Note that $W^\dag_t = \operatorname{Ad}_{\xi_t}^\dag:\mathscr{W}\rightarrow\mathrm{End}(\mathfrak{g})$ satisfies the differential equation
\[ dW^\dag_t = \operatorname{ad}_{X_i}^\dag W^\dag_t \circ db_t^i, \text{ with } W^\dag_0=I, \]
which is linear with smooth coefficients.  Similarly, one may show that 
\[ \overline{W}^\dag := \int_0^\cdot \operatorname{Ad}_{\xi_s}^\dag\,ds \in \mathcal{D}^\infty(\mathscr{H}(\mathrm{End}(\mathfrak{g}))). \]
Also, Theorem \ref{t:inv.smooth} implies that 
\[ \bar{\sigma}_t^{-1} = \left( \int_0^t \operatorname{Ad}_{\xi_s} P\,\operatorname{Ad}_{\xi_s}^\dagger \,ds \right)^{-1} \]
exists and is in $L^{\infty-}(\mu)$ componentwise.  Thus, Equation (\ref{e:Xs})
implies that $\mathbf{X}\in\mathcal{D}^\infty(\mathscr{H})$.

For $f\in C^\infty(G)$ and $(h^1,\ldots,h^k)\in \mathscr{H}$, by Equation \eqref{e:dhf}, 
\[ \partial_h [f(\xi_t)] = ( D[f(\xi_t)],h )_\mathscr{H} 
   = \left\langle \hat{\nabla}f(\xi_t), \int_0^t \operatorname{Ad}_{\xi_s} X_i \dot{h}^i_s\,ds \right\rangle,  \]
and so
\[ \begin{split}
\mathbf{X} [f(\xi_t)] &= ( D[f(\xi_t)],\mathbf{X} )_\mathscr{H} \\
   &= \left\langle \hat{\nabla}f(\xi_t), \int_0^t \operatorname{Ad}_{\xi_s} X_i \left\langle \operatorname{Ad}_{\xi_s}^\dagger \left( \int_0^t \operatorname{Ad}_{\xi_r}
     P\,\operatorname{Ad}_{\xi_r}^\dagger \,dr \right)^{-1} \operatorname{Ad}_{\xi_t}X, X_i\right\rangle ds \right\rangle \\
   &= \left\langle \hat{\nabla}f(\xi_t),  \int_0^t \operatorname{Ad}_{\xi_s} P\operatorname{Ad}^\dagger_{\xi_s}
     \left( \int_0^t \operatorname{Ad}_{\xi_r} P\operatorname{Ad}^\dagger_{\xi_r}\,dr\right)^{-1}
     \operatorname{Ad}_{\xi_t} X \, ds\right\rangle \\
   &= \left\langle \hat{\nabla}f(\xi_t), \operatorname{Ad}_{\xi_t} X \right\rangle
    = \left\langle \tilde{\nabla}f(\xi_t), \operatorname{Ad}_{\xi_t^{-1}}\operatorname{Ad}_{\xi_t} X \right\rangle \\
   &= (\tilde{X} f) (\xi_t),
\end{split} \]
where the penultimate equality used Equation (\ref{e:d2.0}).
\end{proof}

\begin{defn}
For a vector field $\mathbf{X}$ acting on functions of $\mathscr{W}$, denote the
adjoint of $\mathbf{X}$ in the $L^2(\mu)$ inner product by $\mathbf{X}^*$, which has 
domain in $L^2(\mu)$ consisting of functions $G$ such that for all
$F\in\mathcal{D}^{1,2}$,
\[ \mathbb{E}[(\mathbf{X} F)G] \le c\|F\|_{L^2(\mu)} \]
for some constant $c$.  For functions $G$ in the domain of $\mathbf{X}^*$,
\begin{equation*} 
\mathbb{E}[F(\mathbf{X}^*G)] = \mathbb{E}[(\mathbf{X} F)G]
\end{equation*}
for all $F\in\mathcal{D}^{1,2}$.
\end{defn}

Note that for any lifted vector field $\mathbf{X}$ acting on function $F\in\mathcal{D}^{1,2}$ as defined in Definition \ref{d:lift},
\[ \mathbb{E}[\mathbf{X} F] = \mathbb{E}[(DF,\mathbf{X})_\mathscr{H}] = \mathbb{E}[FD^*\mathbf{X}]. \]
Thus, $\mathbf{X}^*=\mathbf{X}^*1=D^*\mathbf{X}$ a.s.  Recall that $D^*$ 
is a continuous operator from $\mathcal{D}^\infty(\mathscr{H})$ into $\mathcal{D}^\infty$;  see for example Theorem V-8.1 and its corollary in \cite{IkedaWatanabe89}.  Thus, 
for $\mathbf{X}$ a vector field on $W$ as defined in Equation 
\eqref{e:lift}, Proposition \ref{p:A.smooth} implies that 
$D^*\mathbf{X}\in\mathcal{D}^\infty$.  This proves the following proposition.

\begin{prop} 
\label{p:K.finite}
Let $\tilde{X}$ be a left invariant vector field on $G$.  Then for the 
vector field on $\mathscr{W}$ defined by
\[ \mathbf{X} = \xi_t'(\omega)^*[\xi_t'(\omega)\xi_t'(\omega)^*]^{-1} \tilde{X}(\xi_t(\omega)), \]
$\mathbf{X}^*\in\mathcal{D}^\infty$, where $\mathbf{X}^*$ is the $L^2(\mu)$-adjoint of $\mathbf{X}$.
\end{prop}

\section{Lie group inequalities}
\label{s:Lie}

Again let $G$ be a Lie group with identity $e$ and Lie algebra $\mathrm{Lie}(G)=\mathfrak{g}$, and suppose $\{X_{i}\}_{i=1}^{k}\subset\mathfrak{g}$ is a H\"ormander set, in the sense of Equation (\ref{e:gen}).  The gradient $\nabla=(\tilde{X}_1,\ldots, \tilde{X}_k)$ and the subLaplacian $L=\sum_{i=1}^k\tilde{X}_i^2$ are operators on smooth functions of $G$ with compact support.  Let $L$ also denote the self-adjoint extension of the subLaplacian and $P_t=e^{tL/2}$ be the heat semigroup as in Definition \ref{d:Pt}.

The following lemmas were proved in \cite{Melcher1} in the context of the Heisenberg Lie group (Lemmas 2.3 and 2.4).  The proofs are identical in the general Lie group case.

\begin{lem} 
\label{l:G.inv}
By the left invariance of $\nabla$ and $P_{t}$, the inequality (\ref{e:MAIN}) 
holds for all $g\in G$, $f\in C_{c}^{\infty}(G)$, and $t>0$, if and only if,
\begin{equation*} 
|\nabla P_{t}f|^{p}(e)\leq K_{p}(t)P_{t}|\nabla f|^{p}(e),
\end{equation*}
for all $f\in C_c^{\infty}(G)$ and $t>0$, where $e\in G$ is the identity element.
\end{lem}

\begin{lem} 
\label{l:rtlt}
For $X\in\mathfrak{g}$,
\[ \tilde{X} P_t f (e) = P_t\hat{X} f (e). \]
for all $f\in C_c^\infty(G)$.  More generally,
\[ \hat{X} P_t f = P_t\hat{X} f, \]
from which the previous equation follows, since $\hat{X}=\tilde{X}$ at $e$.
\end{lem}
 
(The proof of Lemma \ref{l:rtlt}  is actually easier than its analogue Lemma 2.4 in \cite{Melcher1}, since working with functions with compact support -- versus functions with polynomial growth -- requires only the invariance of Haar measure to justify passing the derivative through the integral.)

\subsection{$L^p$-type gradient estimate $(p>1)$} 
\label{s:Lp}

\begin{notation}
For each $r\in\{0,1,\ldots,m\}$, let $\Lambda^r=\Lambda^{k,r}$ be the set of multi-indices $\alpha=(\alpha_0,\alpha_1,\ldots,\alpha_r)\in\{1,\ldots,k\}^{r+1}$.  For any $\alpha\in\Lambda^{r}$, define
\[  \alpha' := (\alpha_1,\ldots,\alpha_r) \text { and }\]
\[  \overline{\alpha} := (\alpha_r,\ldots,\alpha_0) = \alpha \text{ reversed }. \]
Define the order of $\alpha$ by $|\alpha|:=r+1$.  Let
\[ X_\alpha = [X_{\alpha_r},[\cdots,[X_{\alpha_1},X_{\alpha_0}]\cdots]] 
   = \operatorname{ad}_{X_{\alpha_r}}\cdots \operatorname{ad}_{X_{\alpha_1}} X_{\alpha_0} \text{ and }\]
\[ X^\alpha=X_{\alpha_r}\cdots X_{\alpha_0}. \] 
When $r=0$ and $|\alpha|=1$, that is, $\alpha=(\alpha_0)$, then $X^\alpha=X_{\alpha_0}=X_\alpha$.  For each $\alpha\in\Lambda^{r}$, there exist $\epsilon_{\beta,\alpha}\in\mathbb{Z}$ such that
\[ X_\alpha = \sum_{\beta\in\Lambda^{r}} \epsilon_{\beta,\alpha} X^\beta. \] 
\end{notation}

\begin{prop}
\label{r:K.finite}
For any $X\in\mathfrak{g}$, $\hat{X}$ may be written as
\begin{equation}
\label{e:rt.lt}
\hat{X} = \sum_{r=0}^m \sum_{\alpha\in\Lambda^{r}} c_{\alpha}\tilde{X}^\alpha,
\end{equation}
with $c_{\alpha}:G\rightarrow \mathbb{R}$ (some of these are 0) such that  $c_{\alpha}(\xi_t)\in \mathcal{D}^\infty$, for all $t\in[0,1]$.  
\end{prop}

\begin{proof}
Recall from Notation \ref{n:Sigma} that 
\[ \begin{split}
 \Sigma_r &= \{[X_{i_1},[\cdots,[X_{i_{r-1}},X_{i_r}]\cdots]: i_1,\ldots,i_r     \in\{1,\ldots,k\}\} \\
   &= \{X_\alpha : \alpha\in\Lambda^{r} \}, 
\end{split} \]
for $r=0,\ldots,m$.  Recall also from Notation \ref{n:Sigma} that $\{X_i,Y_j : i\in\{1,\ldots,k\}, j\in\{1,\ldots,d-k\}\}$ of $\mathfrak{g}$ is an orthonormal basis, where $d=\mathrm{dim}(G)$ and, for each $j\in\{1,\ldots,d-k\}$, $Y_j$ is some commutator $X_{\alpha(j)}\in\Sigma_{r(j)}$ for some $\alpha(j)\in\Lambda^{r(j)}$, $r(j)\in\{1,\ldots,m\}$.  Thus, for any $g\in G$ and $X\in\mathfrak{g}$,
\begin{align*}
\hat{X}(g) &= R_{g*}X = L_{g*}L_{g^{-1}*}R_{g*}X = L_{g*}\operatorname{Ad}_{g^{-1}}X \\
   &= L_{g*} \left( \sum_{i=1}^k \left\langle \operatorname{Ad}_{g^{-1}}X,X_i\right\rangle X_i 
      + \sum_{j=1}^{d-k} \left\langle\operatorname{Ad}_{g^{-1}}X,Y_j\right\rangle Y_j \right) \\
   &=  L_{g*} \left( \sum_{i=1}^k \left\langle \operatorname{Ad}_{g^{-1}}X,X_i\right\rangle X_i 
      + \sum_{j=1}^{d-k} \sum_{\alpha\in\Lambda^{r(j)} } \epsilon_{\alpha,\alpha(j)} 
       \left\langle\operatorname{Ad}_{g^{-1}}X,Y_j\right\rangle X^\alpha \right) \\
   &=  \sum_{i=1}^k \left\langle \operatorname{Ad}_{g^{-1}}X,X_i\right\rangle \tilde{X}_i(g) 
      + \sum_{j=1}^{d-k} \sum_{\alpha\in\Lambda^{r(j)} } \epsilon_{\alpha,\alpha(j)} 
       \left\langle\operatorname{Ad}_{g^{-1}}X,Y_j\right\rangle \tilde{X}^\alpha(g)
\end{align*}
where $\epsilon_{\alpha,\alpha(j)}\in\mathbb{Z}$.  So 
\[ \hat{X}(g) = \sum_{r=0}^m \sum_{\alpha\in\Lambda^{r}} c_{\alpha}\tilde{X}^\alpha(g), \]
where 
\[ c_{\alpha}(g) = \left\{\begin{array}{ll}
	\left\langle \operatorname{Ad}_{g^{-1}}X,X_i\right\rangle & \text{ when } r=0 \text{ and } \alpha=(i) \\
	\epsilon \left\langle\operatorname{Ad}_{g^{-1}}X,Y_j\right\rangle, \epsilon\in\mathbb{Z} & \text{ when } r\in\{1,\ldots,m\} 
	\end{array}\right. .\]

Note that $\operatorname{Ad}_{\xi_t}$ satisfies the Stratonovich stochastic differential equation
\begin{equation*}
d\operatorname{Ad}_{\xi} = \operatorname{Ad}_{\xi} \circ\operatorname{ad}_{db} = \operatorname{Ad}_{\xi} \operatorname{ad}_{X_i} \circ db^i, ~ \text{ with } \operatorname{Ad}_{\xi_0}=I.
\end{equation*}
By differentiating the identity $\operatorname{Ad}_{\xi_t}\operatorname{Ad}_{\xi_t}^{-1} = I$, one may verify that $\operatorname{Ad}_{\xi_t}^{-1}=\operatorname{Ad}_{\xi_t^{-1}}$ satisfies 
\begin{equation*}
d\operatorname{Ad}_{\xi^{-1}} = -\circ\operatorname{ad}_{db} \operatorname{Ad}_{\xi^{-1}}  
   = -  \operatorname{ad}_{X_i} \operatorname{Ad}_{\xi^{-1}} \circ db^i, ~  \text{ with } \operatorname{Ad}_{\xi^{-1}_0}=I
\end{equation*}
a linear differential equation with smooth coefficients.  Then by Theorem V-10.1 of Ikeda and Watanabe \cite{IkedaWatanabe89}, $\operatorname{Ad}_{\xi_t^{-1}}\in\mathcal{D}^\infty(\mathrm{End}(\mathfrak{g}))$ componentwise 
with respect to some basis. 

The function $u:\mathrm{End}(\mathfrak{g})\rightarrow\mathbb{R}$ given by $u(W) = \left\langle W X,Y\right\rangle$ is a smooth function for any fixed $X,Y\in\mathfrak{g}$.  Thus, $u(\operatorname{Ad}_{\xi^{-1}_t})\in\mathcal{D}^\infty$ for all $t\in[0,1]$.  Since $c_{\alpha}(\xi_t)=\epsilon u(\operatorname{Ad}_{\xi^{-1}_t})$, with $Y=X_i$ or $Y_j$, this implies that $c_{\alpha}(\xi_t)\in\mathcal{D}^\infty$, for all $\alpha\in\Lambda^{r}$.
\end{proof}

\begin{thm}
\label{t:Lie}
For all $p\in(1,\infty)$, $K_p(t)<\infty$, where $K_p(t)$ are the functions defined in Notation \ref{n:best}.
\end{thm}

\begin{proof}
Lemma \ref{l:G.inv} implies that the inequality (\ref{e:MAIN}) is translation invariant on
groups.  Thus it suffices to determine a finite coefficient $K_p(t)$ such that the inequality holds at the identity.  

Note that for any $X\in\mathfrak{g}$, Lemma \ref{l:rtlt} and Equation (\ref{e:rt.lt}) imply that
\[ |\tilde{X} P_t f|^2(e) =  |\hat{X} P_t f|^2(e) =  |P_t \hat{X} f|^2(e) 
   \le C \sum_{r=0}^m \sum_{\alpha\in\Lambda^{r}} 
     |P_t c_{\alpha} \tilde{X}^\alpha f|^2(e), \]
for a constant $C=C(k,m)$.  Equation (\ref{e:integral}) implies that, for any $f\in C_c^\infty(G)$, $P_tf(e)=\mathbb{E}[f(\xi_t)]$, where $\xi$ is the solution to the Stratonovich equation (\ref{e:SDE}).  Thus, for any $\alpha\in\Lambda^{r}$, 
\[ 
\begin{split}
|P_t c_\alpha\tilde{X}^\alpha f|(e)
   &\le \mathbb{E}|c_\alpha(\xi_t)(\tilde{X}^\alpha f)(\xi_t)| 
   = \mathbb{E}|c_\alpha(\xi_t) \mathbf{X}^{\alpha'}[(\tilde{X}_{\alpha_0}f)(\xi_t)]| \\
   &= \mathbb{E}\left|\left(\mathbf{X}^{\overline{\alpha'}}\right)^*[c_\alpha(\xi_t)] (\tilde{X}_{\alpha_0}f)(\xi_t)\right| \\
   &\le \left( \mathbb{E}\left|\left(\mathbf{X}^{\overline{\alpha'}}\right)^*[c_\alpha(\xi_t)]\right|^q\right)^{1/q}
     \left( \mathbb{E}|(\tilde{X}_{\alpha_0}f)(\xi_t)|^p \right)^{1/p} \\
   &= \left( \mathbb{E}\left|\left(\mathbf{X}^{\overline{\alpha'}}\right)^*[c_\alpha(\xi_t)]\right|^q\right)^{1/q}
     \left( P_t|\tilde{X}_{\alpha_0}f|^p(e) \right)^{1/p} \\
   &\le \left( \mathbb{E}\left|\left(\mathbf{X}^{\overline{\alpha'}}\right)^*[c_\alpha(\xi_t)]\right|^q\right)^{1/q}
     \left( P_t|\nabla f|^p(e) \right)^{1/p} ,
\end{split} 
\]
by H\"older's inequality, where $q$ is the conjugate exponent to $p$, $\mathbf{X}^\alpha$ is the lifted vector field on $W$ of the vector field $\tilde{X}^\alpha$, as defined in Equation (\ref{e:lift}), and 
$(\mathbf{X}^\alpha)^*=\mathbf{X}_{\alpha_r}^*\cdots\mathbf{X}_{\alpha_0}^*$
(so $\left(\mathbf{X}^{\overline{\alpha'}}\right)^* = \mathbf{X}_{\alpha_1}^*\cdots\mathbf{X}_{\alpha_r}^*$).
Propositions \ref{p:K.finite} and \ref{r:K.finite} imply that $\left(\mathbf{X}^{\overline{\alpha'}}\right)^*[c_{\alpha}(\xi_t)]\in L^{\infty-}(\mu)$, for all $\alpha\in\Lambda^{r}$.

So in particular, using the above with $X=X_i$ gives
\[  \begin{split}
|\nabla P_t f|^p (e)
   &= \left(\sum_{i=1}^k |\tilde{X}_i P_t f|^2(e)  \right)^{p/2} \\
   &\le C\left[\sum_{i=1}^k \sum_{r=0}^m \sum_{\alpha\in\Lambda^{r}}
      \left( \mathbb{E}\left|\left(\mathbf{X}^{\overline{\alpha'}}\right)^*\left[c_{i,\alpha}(\xi_t)\right]\right|^q\right)^{p/q}\right]
       P_t|\nabla f|^p (e), 
\end{split} 
\]
where $C=C(k,m,p)$ and $q=\frac{p}{p-1}$.  Thus, the inequality (\ref{e:MAIN}) holds with 
\begin{equation} 
\label{e:K_p}
C_p(t) = C(k,m,p)\sum_{i=1}^k \sum_{r=0}^m \sum_{\alpha\in\Lambda^{r}}
      \left( \mathbb{E}\left|\left(\mathbf{X}^{\overline{\alpha'}}\right)^*\left[c_{i,\alpha}(\xi_t)\right]\right|^q\right)^{p/q}. 
\end{equation}
Therefore, $K_p(t)\le C_p(t)<\infty$ for all $t>0$ and $p\in(1,\infty)$.
\end{proof}

It is important to note that, in this general Lie group case, there is currently no good control over the behavior of the functions $C_p$ in Equation (\ref{e:K_p}) with respect to $t$.  In fact, from certain scaling arguments, it is expected that $C_p(t)\rightarrow\infty$ as $t\rightarrow0$; see for example \cite{Bell5,KS2}.  However, these coefficients are almost certainly not optimal.

To explore cases where the behavior of these coefficients is more understood, it will become useful to extend the set of test functions considered.  The following proposition relaxes the condition of compact support to boundedness with bounded first order derivatives.

\begin{prop}
\label{p:bounded}
For all $p\in(1,\infty)$,
\begin{equation*}
|\nabla P_t f|^p \le K_p(t) P_t |\nabla f |^p,
\end{equation*}
for all $f\in C_b^\infty(G)$ with bounded derivatives of first order and $t>0$.
\end{prop}

\begin{proof}
Let $f\in C_b^\infty(G)$ with bounded first order derivatives, and let $\varphi_m\in C_c^\infty(G,[0,1])$ be a sequence of functions such that $\varphi_m\uparrow1$, $\varphi_m(g)=1$ when $|g|\le m$ (for some norm on $G$), and $\sup_m\sup_{g\in G} |\tilde{X} \varphi_m|<\infty$ for all $X\in\mathfrak{g}$; see Lemma 3.6 of \cite{Driver92b}.  Then $f_m=\varphi_m f\in C_c^\infty(G)$, and so there exists an optimal function $K_p(t)<\infty$ such that
\[ |\nabla P_t f_m|^p \le K_p(t) P_t |\nabla f_m |^p. \]
for all $t>0$.  For any $X\in\mathfrak{g}$,
\[ \begin{split}
\lim_{m\rightarrow\infty} |\tilde{X} f_m - \tilde{X} f| 
   &= \lim_{m\rightarrow\infty} |(\tilde{X}\varphi_m) f + \varphi_m\tilde{X} f - \tilde{X} f|  \\
   &\le \lim_{m\rightarrow\infty} |\tilde{X} \varphi_m| |f| + |\varphi_m-1||\tilde{X} f|
   = 0
\end{split} \]
implies that $|\nabla f_m|\rightarrow|\nabla f|$ boundedly.  
Thus, by the dominated convergence theorem,
\[  \lim_{m\rightarrow\infty} P_t |\nabla f_m|^p = P_t |\nabla f|^p. \]
Similarly,
\[\begin{split}
\lim_{m\rightarrow\infty} |\tilde{X}P_tf_m - \tilde{X}P_tf| 
   &= \lim_{m\rightarrow\infty} |P_t\hat{X}f_m - P_t\hat{X}f| \\
   &\le  \lim_{m\rightarrow\infty} P_t|\hat{X}f_m - \hat{X}f| \\
   &\le  \lim_{m\rightarrow\infty} P_t(|\hat{X}\varphi_m||f|) + P_t(|\varphi_m-1||\hat{X}f|) 
   = 0
\end{split}\]
by dominated convergence, and hence
\[ \lim_{m\rightarrow\infty} |\nabla P_t f_m| = |\nabla P_t f|. \]
Thus,
\[\begin{split}
|\nabla P_t f|^p 
   = \lim_{m\rightarrow\infty} |\nabla P_tf_m|^p 
   \le K_p(t) \lim_{m\rightarrow\infty} P_t |\nabla f_m|^p 
   =  K_p(t) P_t |\nabla f|^p.
\end{split}\]
\end{proof}

\subsection{Poincar\'e inequality}
\label{s:poincare}

The following result is a direct corollary to Theorem \ref{t:Lie}.  The proof is completely analogous to the proof of Theorem 4.2 in \cite{Melcher1} in the Heisenberg Lie group context.

\begin{thm}[Poincar\'e Inequality]
\label{t:Poincare}
Let $K_2(t)$ be the best function for which (\ref{e:MAIN}) holds for $p=2$, and let $p_t(g)\,dg$ be the hypoelliptic heat kernel.  Then 
\begin{equation}
\label{e:Poincare}
\int_G f^2(g) p_t(g)\, dg - \left( \int_G f(g)p_t(g)\,dg \right)^2 \le \Lambda(t)\int_G |\nabla f|^2(g)p_t(g)\,dg,
\end{equation}
for all $f\in C_c^\infty(G)$ and $t>0$, where
\[ \Lambda(t) = \int_0^t K_2(s)\,ds. \]
\end{thm}
\begin{proof}
Let $F_t(g) = (P_tf)(g)$.  Then
\[
\frac{d}{ds} P_{t-s}F_s^2 = P_{t-s} \left( -\frac{1}{2}LF_s^2 + F_sLF_s \right) = -P_{t-s}|\nabla F_s|^2.
\]
Integrating this equation on $s$ implies that
\begin{equation*}
\begin{split}
P_t f^2 - (P_tf)^2 &=  \int_0^t P_{t-s}|\nabla F_s|^2 \,ds \\
   &=  \int_0^t P_{t-s}|\nabla P_sf|^2 \,ds \\
   &\le \int_0^t K_2(s) P_{t-s}P_s|\nabla f|^2 \,ds
    = \left(\int_0^t K_2(s)\, ds \right) \cdot P_t|\nabla f|^2
\end{split}
\end{equation*}
where the inequality follows from Theorem \ref{t:Lie}.  Evaluating the above at $e\in G$ gives the desired result.
\end{proof}

This theorem is less useful in the general Lie group case because nothing is known about the integrability of $K_p(t)$.  However, the next two sections show that, when $G$ is a nilpotent Lie group, $K_p(t)$ is a bounded function for all $p\in(1,\infty)$.  In particular, when $p=2$, this implies the Poincar\'{e} inequality holds with $\Lambda(t)<\infty$, for all $t>0$.

\subsubsection{Stratified nilpotent Lie groups}
\label{ss:stratified}

\begin{defn}
A Lie algebra $\mathfrak{g}$ is said to be {\em nilpotent} if $\operatorname{ad}_X$ is a nilpotent endomorphism of $\mathfrak{g}$ for all $X\in\mathfrak{g}$, that is, if there exists $m\in\mathbb{N}$ such that
\[ 
\operatorname{ad}_{Y_1}\cdots\operatorname{ad}_{Y_{m-1}}Y_m = [Y_1,[\cdots,[Y_{m-1},Y_m]\cdots] = 0, 
\]
for any $Y_1,\ldots,Y_m\in\mathfrak{g}$.  If $m$ is the smallest number for which the above equality holds, $\mathfrak{g}$ is nilpotent of step $m$.  A Lie group $G$ is nilpotent if $\mathfrak{g}=\mathrm{Lie}(G)$ is a nilpotent Lie algebra.  
\end{defn}

\begin{defn} 
\label{d:dilation}
A family of {\em dilations} on a Lie algebra $\mathfrak{g}$ is 
a family of algebra automorphisms $\{\Phi_r\}_{r>0}$ on $\mathfrak{g}$ of 
the form $\Phi_r = \exp(W\log r),$ where $W$ is a diagonalizable 
linear operator on $\mathfrak{g}$ with positive eigenvalues.
\end{defn}


\begin{defn}
A {\em stratified} group $G$ is a simply connected nilpotent group for which there exists a subset of the Lie algebra $V_1\subset\mathfrak{g}$, such that $\mathfrak{g} = \oplus_{j=1}^m V_j$ with $V_{j+1}=[V_1,V_j]$, for $j=1,\ldots,m-1$, and $V_{m+1}=[V_1,V_m]=\{0\}$.
\end{defn}

For a general exposition on nilpotent Lie groups and dilations, see \cite{FollandStein82,Goodman} and references contained therein.  If $G$ is a stratified Lie group, a natural family of dilations may be defined on $\mathfrak{g}$ by setting $\Phi_r(X)=r^jX$, for all $X\in V_j$.  The generator $W$ of this dilation acts on parts of the vector space decomposition by $WV_j = jV_j$, for each $j=1,\ldots,m$.  The automorphism $\Phi_r$ induces a group dilation $\phi_r$ via the exponential maps, $\phi_r=\exp\circ\Phi_r\circ\exp^{-1}$.  Since $G$ is a simply connected nilpotent group, the exponential map is in fact a global diffeomorphism on $\mathfrak{g}$, and $\exp^{-1}$ exists everywhere on $G$; see for example Theorem 3.6.2 of Varadarajan \cite{Varadarajan}.  Then for each $X\in V_1$,
\begin{multline} 
\label{e:dilation}
\tilde{X}(f\circ \phi_r)(g) = \frac{d}{d\epsilon}\bigg|_0 (f\circ\phi_r)(ge^{\epsilon X}) 
   = \frac{d}{d\epsilon}\bigg|_0 f(\phi_r(g)\phi_r(e^{\epsilon X})) \\
   = \frac{d}{d\epsilon}\bigg|_0 f(\phi_r(g)e^{r\epsilon X}) 
   = \frac{d}{d\epsilon}\bigg|_0 rf(\phi_r(g)e^{\epsilon X}) 
   = r (\tilde{X} f\circ \phi_r)(g),
\end{multline}
for all $f\in C^1(G)$, where the second equality used that$\phi_r$ is a homomorphism.  Let $\{X_i\}_{i=1}^k\subset V_1$ be a basis of $V_1$, and consider the operators $\nabla=(\tilde{X}_1,\ldots,\tilde{X}_k)$ and $L=\sum_{i=1}^k \tilde{X}_i^2$.  Equation (\ref{e:dilation}) implies that 
\begin{equation}
\label{e:gradscale} 
\nabla(f\circ\phi_r) = r(\nabla f)\circ\phi_r, 
\end{equation}
and thus the following proposition.
\begin{prop}
Let $L$ denote the self-adjoint extension of $\sum_{i=1}^k\tilde{X}_i^2$, and $P_t=e^{tL/2}$ be as in Definition \ref{d:Pt}.  Then 
\[ L(f\circ\phi_r) = r^2(Lf)\circ\phi_r \]
and
\begin{equation}
\label{e:Pscale}
P_t(f\circ\phi_r) = e^{tL/2}(f\circ\phi_r) = \left( e^{r^2tL/2}f \right) \circ \phi_r 
	= \left( P_{r^2t}f \right)\circ\phi_r,
\end{equation}
for any $f\in C_c^\infty(G)$.
\end{prop}

\begin{proof}
Let $\mathcal{E}^0(f,h):=\sum_{i=1}^k (\tilde{X}_i f,\tilde{X}_i h)_{L^2(G)}$ be a Dirichlet form associated to $L$.  Recall from Section \ref{s:intro} that $\mathcal{E}^0$ has a closed extension $\mathcal{E}$.  By definition,
\[ f_1\in C_c^\infty(G) \text{ and } Lf_1=h \iff \mathcal{E}(f_1,f_2) = (h,f_2), ~\forall f_2\in\mathrm{Dom}(\mathcal{E}). \]
Now note that
\begin{align*}
\mathcal{E}^0(f\circ\phi_r,f\circ\phi_r) 
   &= \sum_{i=1}^k \int_G |\tilde{X}_i(f\circ\phi_r)|^2(g)\,dg \\
   &= \sum_{i=1}^k r^2\int |(\tilde{X}_i f)\circ\phi_r|^2(g)\,dg \\
   &= \sum_{i=1}^k r^2\int |\tilde{X}_i f|^2(g) J(r^{-1}) \,dg = r^2J(r^{-1})\mathcal{E}^0(f,f),
\end{align*}
where $J(r)$ is the Jacobian of the transformation $\phi_r$, 
\[ J(r) = \prod_{j=1}^m (r^j)^{d_j} \]
with $d_j=\dim(V_j)$.  Thus, $J(r^{-1})=J(r)^{-1}$.  So $f\in\mathrm{Dom}(\mathcal{E})$ implies that $f\circ\phi_r\in\mathrm{Dom}(\mathcal{E})$, and, in general, $\mathcal{E}(f\circ\phi_r,h\circ\phi_r) = r^2J(r^{-1})\mathcal{E}(f,h)$, for $f,h\in\mathrm{Dom}(\mathcal{E})$.  Replacing $h$ here by $h\circ\phi_{r^{-1}}$ gives
\begin{align*} 
\mathcal{E}(f\circ\phi_r,h) &= r^2J(r^{-1})\mathcal{E}(f,h\circ\phi_{r^{-1}}) \\
   &= r^2J(r^{-1})(Lf,h\circ\phi_{r^{-1}})_{L^2(G)} \\
   &= r^2J(r^{-1})J(r)(Lf\circ\phi_r,h)_{L^2(G)} 
   = r^2(Lf\circ\phi_r,h)_{L^2(G)},
\end{align*}
implies that if $f\in\mathrm{Dom}(L)$, then $f\circ\phi_r\in\mathrm{Dom}(L)$ and $L(f\circ\phi_r) = r^2Lf\circ\phi_r$.  

Now, for $r>0$, let $U_r:L^2(G)\rightarrow L^2(G)$ be the unitary operator given by $U_rf=\frac{1}{\sqrt{J(r^{-1})}}f\circ\phi_r$.  Then
\[ LU_r = r^2U_r L = U_r(r^2L) \]
as operators, and thus $U_r^{-1}LU_r=r^2L$.  Then 
\[ U_r^{-1}e^{tL/2} U_r  = e^{tU_r^{-1}LU_r/2}= e^{r^2tL/2}, \]
from which it follows that
\[ r^2e^{tL/2}(f\circ\phi_r) = e^{tL/2} U_rf = U_r e^{r^2tL/2}f = r^2(e^{r^2tL/2}f)\circ \phi_r. \]
\end{proof} 

This give the following proposition.

\begin{prop} 
\label{p:t.ind}
Suppose $G$ is a stratified Lie group with vector space decomposition $\oplus_{j=1}^m V_j$.  Let $\{X_i\}_{i=1}^k\subset V_1$, $\nabla$, and $L$ be as above, and let $p\in(1,\infty)$.  If 
$K_p$ is the best constant such that
\[ |\nabla P_1f|^p \le K_p P_1|\nabla f|^p, \]
for all $f\in C_c^\infty(G),$  then $K_p(t)=K_p$ for all $t>0$, where
$K_p(t)$ is the function defined in Notation \ref{n:best}.
\end{prop}

\begin{proof}
By Equations (\ref{e:gradscale}) and (\ref{e:Pscale}),
\[
\begin{split}
|\nabla P_{t}(f\circ\phi_{t^{-1/2}})|^{p}  
	&=|\nabla\lbrack(P_{1}f)\circ \phi_{t^{-1/2}}]|^{p}
	= |t^{-1/2}(\nabla P_{1}f)\circ\phi_{t^{-1/2}}|^{p}\\
	&\leq K_{p}t^{-p/2}\left(  P_{1}|\nabla f|^{p}\right)  \circ\phi_{t^{-1/2}}
	= K_{p}t^{-p/2}P_{t}\left(  |\nabla f|^{p}\circ\phi_{t^{-1/2}}\right) \\
	&=K_{p}P_{t}\left(  |\nabla f\circ\phi_{t^{-1/2}})|^{p}\right).
\end{split}
\]
Replacing $f$ by $f\circ\phi_{t^{1/2}}$ in the above computation proves the assertion.
Moreover, reversing the above argument shows that $|\nabla P_{t}f|^{p}\leq
K_{p}P_{t}|\nabla f|^{p}$ implies that $|\nabla P_{1} f|^{p}\leq K_{p}P_{1}|\nabla
f|^{p}.$
\end{proof}

\subsubsection{Nilpotent Lie groups}
\label{ss:nilpotent}

Now let $G$ be a general nilpotent Lie group.  Because not all nilpotent Lie groups admit dilations, the functions $K_p(t)$ are not scale invariant in this context.  However, covering $G$ with a group which has a family of dilations adapted to its structure, shows that there exists some constant $K_p<\infty$ for which $K_p(t)<K_p$ for all $t>0$.

\begin{defn}
Let $\mathcal{L}=\mathcal{L}(k,m)$ be the {\em free nilpotent Lie algebra of step $m$ with $k$ generators} $\{e_i\}_{i=1}^k$.  Then $\mathcal{L}$ is the unique (up to isomorphism) nilpotent Lie algebra of rank $m$ such that, for every nilpotent Lie algebra $\mathfrak{g}$ of rank $m$ and map $\tilde{\Pi}:\{e_1,\ldots,e_k\}\rightarrow\mathfrak{g}$, there exists a unique homomorphism $\Pi:\mathcal{L}\rightarrow\mathfrak{g}$ which extends $\tilde{\Pi}$.  Let $\mathcal{N}=\mathcal{N}(k,m)$ be the {\em free nilpotent Lie group of rank $m$ with $k$ generators}, which is the simply connected group of $\mathcal{L}(k,m)$.  
\end{defn}

The Lie algebra $\mathcal{L}(k,m)$ admits a vector space decomposition by setting $V_1=\mathrm{span}\{e_1,\ldots,e_k\}$.  Thus, $\mathcal{N}$ is a stratified Lie group with H\"ormander set $\{e_i\}_{i=1}^k\subset\mathcal{L}$; for definitions and further details, see \cite{Varopoulos92}.  Let  $\nabla_\mathcal{L}=(\tilde{e}_1,\ldots,\tilde{e}_k)$,  $\mathscr{L}=\sum_{i=1}^k \tilde{e}_i^2$, and $\mathscr{P}_t=e^{t\mathscr{L}/2}$.  Theorem \ref{t:Lie} and Proposition \ref{p:t.ind} imply that,
for all $p\in(1,\infty)$, there exist constants $K_p^\mathcal{L}<\infty$ such that
\begin{equation}
\label{e:free}
|\nabla_\mathcal{L} \mathscr{P}_t f|^p \le K_p^\mathcal{L} \mathscr{P}_t |\nabla_\mathcal{L} f |^p,
\end{equation}
for all  $f\in C_c^\infty(\mathcal{N})$ and $t>0$.

\begin{prop}
\label{p:nilpotent}
Let $G$ be a nilpotent group of step $m$ with H\"ormander set $\{X_i\}_{i=1}^k$.  Then $K_p(t)\le K_p^\mathcal{L}$ for all $t>0$, where $K_p(t)$ is the function defined in Notation \ref{n:best}.
\end{prop}

\begin{proof}
By definition of $\mathcal{L}=\mathcal{L}(k,m)$, there exists a unique Lie algebra homomorphism $\Pi:\mathcal{L}\rightarrow\mathfrak{g}$ such that $\Pi(e_i) = X_i$.  Then $\Pi$ induces a group homomorphism $\pi:\mathcal{N}\rightarrow G$ via the exponential maps,
\[ \pi = \exp_G \circ \Pi \circ \exp_{\mathcal{N}}^{-1}. \]  
Again, because $\mathcal{N}$ is a simply connected nilpotent Lie group, the exponential map on $\mathcal{L}$ is a global diffeomorphism.  Note that $\pi_*=\Pi$,
\begin{equation*}
\begin{CD}
\mathcal{L}(k,m) @>\Pi>> \mathfrak{g} \\
   @V{\exp_\mathcal{N}}VV                        @VV{\exp_G}V \\
   \mathcal{N}(k,m) @>>\pi> G
\end{CD}
\end{equation*}
and the vector fields $\tilde{X}_i$ and $\tilde{e}_i$ are $\pi$-related; that is,
\[  \tilde{e}_\alpha (f \circ \pi)=(\tilde{X}_\alpha f)\circ\pi, \] 
for any multi-index $\alpha\in\Lambda^{r}$ and $f\in C_c^\infty(G)$.  Note that $f\circ\pi\in C_b^\infty(\mathcal{N})$ and has bounded first order derivatives.  Thus, by Proposition \ref{p:bounded},
\[ 
|\nabla P_t f|^p(e) = |\nabla_\mathcal{L} \mathscr{P}_t (f\circ\pi)|^p(e_\mathcal{N})
   \le K_p^\mathcal{L} \mathscr{P}_t |\nabla_\mathcal{L}(f\circ\pi) |^p(e_\mathcal{N}) 
   = K_p^\mathcal{L} P_t|\nabla f|^p(e),
\]
where $e_\mathcal{N}$ is the identity element of $\mathcal{N}$.  Since $K_p(t)$ is the
best constant for which 
\[ |\nabla P_t f|^p (e)\le K_p(t) P_t|\nabla f|^p(e) \] 
holds, the above implies that $K_p(t)\le K_p^\mathcal{L}$ for all $t>0$.
\end{proof}

This method of lifting the vector fields to a free nilpotent Lie algebra was learned from  \cite{Varopoulos86,Varopoulos92}.  A generalization of this procedure may be found in \cite{Stein76}.

\begin{remark}
\label{r:elliptic}
Note that the above argument is independent of the minimality of the H\"ormander set $\{X_i\}_{i=1}^k$.
So suppose that the collection $\{X_i\}_{i=1}^k$ spans the Lie algebra $\mathfrak{g}$.  Since $G$ is a nilpotent Lie group (and thus unimodular) it is then well known that the 
operator $L=\sum_{i=1}^k \tilde{X}_i^2$ is in fact the Laplace-Beltrami operator on the Riemannian manifold $(G,\left\langle\cdot,\cdot\right\rangle)$.  Then it is well known that the inequality (\ref{e:MAIN}) holds with exponential coefficients:
\[ |\nabla P_t f|^p \le e^{pkt} P_t|\nabla f|^p, \]
where $-2k$ is a lower bound on the Ricci curvature; see for example Theorem 1.1 in \cite{Melcher1}.  Proposition \ref{p:nilpotent} improves this result by implying that there exists a $K_p<\infty$ independent of $t$ such that
\[ |\nabla P_t f|^p \le K_pP_t|\nabla f|^p, \]
for all $f\in C_p^\infty(G)$ and $t>0$.  This implies the following corollary.
\end{remark}

\begin{cor}
\label{c:elliptic}
Let $G$ be a nilpotent Lie group of step $m$ and $\{X_i\}_{i=1}^k\subset\mathfrak{g}$ such that $\{X_i\}_{i=1}^k$ spans the Lie algebra $\mathfrak{g}$.  Then, for $K_p(t)$ as in Notation \ref{n:best},
\[ K_p(t) \le \min\{K_p^\mathcal{L},e^{pkt}\}, \]
where $K_p^{\mathcal{L}}$ is the best constant so that (\ref{e:MAIN}) holds on $\mathcal{L}(k,m)$ and $-2k$ is a lower bound on the Ricci curvature associated to the Riemannian metric determined by $L= \sum_{i=1}^k \tilde{X}_i^2$.
\end{cor}

This also gives the following Poincar\'{e} inequality for nilpotent Lie groups.

\begin{cor}
\label{c:nilPoincare} 
Suppose $G$ is a nilpotent Lie group, and let $K_2$ be a finite constant for which
the inequality (\ref{e:MAIN}) holds for $p=2$.  Then the inequality (\ref{e:Poincare})
holds with $\Lambda(t)=K_2t$, for all $t>0$.
\end{cor}

\bibliographystyle{amsplain}
\bibliography{biblio}

\providecommand{\bysame}{\leavevmode\hbox to3em{\hrulefill}\thinspace}
\providecommand{\MR}{\relax\ifhmode\unskip\space\fi MR }
\providecommand{\MRhref}[2]{%
  \href{http://www.ams.org/mathscinet-getitem?mr=#1}{#2}
}
\providecommand{\href}[2]{#2}
\begin{thebibliography}{10}

\bibitem{Coulhon}
Pascal Auscher, Thierry Coulhon, Xuan~Thinh Duong, and Steve Hofmann,
  \emph{Riesz transform on manifolds and heat kernel regularity}, Ann. Sci.
  \'Ecole Norm. Sup. (4) \textbf{37} (2004), no.~6, 911--957. \MR{MR2119242
  (2005k:58043)}

\bibitem{Bakry90}
Dominique Bakry, \emph{Ricci curvature and dimension for diffusion semigroups},
  Stochastic processes and their applications in mathematics and physics
  (Bielefeld, 1985), Math. Appl., vol.~61, Kluwer Acad. Publ., Dordrecht, 1990,
  pp.~21--31. \MR{92e:58231}

\bibitem{BakryEmery84}
Dominique Bakry and Michel {\'E}mery, \emph{Hypercontractivit\'e de
  semi-groupes de diffusion}, C. R. Acad. Sci. Paris S\'er. I Math.
  \textbf{299} (1984), no.~15, 775--778. \MR{86f:60097}

\bibitem{BakryEmery85}
\bysame, \emph{Diffusions hypercontractives}, S\'eminaire de probabilit\'es,
  XIX, 1983/84, Lecture Notes in Math., vol. 1123, Springer, Berlin, 1985,
  pp.~177--206. \MR{88j:60131}

\bibitem{Bell5}
Denis~R. Bell and Salah Eldin~A. Mohammed, \emph{The {M}alliavin calculus and
  stochastic delay equations}, J. Funct. Anal. \textbf{99} (1991), no.~1,
  75--99. \MR{92k:60124}

\bibitem{Coulhon03}
Thierry Coulhon and Xuan~Thinh Duong, \emph{Riesz transform and related
  inequalities on noncompact {R}iemannian manifolds}, Comm. Pure Appl. Math.
  \textbf{56} (2003), no.~12, 1728--1751. \MR{2 001 444}

\bibitem{Davies}
Edward~Brian Davies, \emph{One-parameter semigroups}, London Mathematical
  Society Monographs, vol.~15, Academic Press Inc. [Harcourt Brace Jovanovich
  Publishers], London, 1980. \MR{82i:47060}

\bibitem{Driver04}
Bruce~K. Driver, \emph{Curved {W}iener space analysis}, Real and stochastic
  analysis, Trends Math., Birkh\"auser Boston, Boston, MA, 2004, pp.~43--198.
  \MR{MR2090752 (2005g:58066)}

\bibitem{Driver92b}
Bruce~K. Driver and Leonard Gross, \emph{Hilbert spaces of holomorphic
  functions on complex {L}ie groups}, New trends in stochastic analysis
  (Charingworth, 1994), World Sci. Publishing, River Edge, NJ, 1997,
  pp.~76--106. \MR{2000h:46029}

\bibitem{Melcher1}
Bruce~K. Driver and Tai Melcher, \emph{Hypoelliptic heat kernel inequalities on
  the {H}eisenberg group}, Journal of Functional Analysis \textbf{221} (2005),
  no.~2, 340--365.

\bibitem{FollandStein82}
Gerald~B. Folland and Elias~M. Stein, \emph{Hardy spaces on homogeneous
  groups}, Princeton University Press, Princeton, N.J., 1982.

\bibitem{Goodman}
Roe~W. Goodman, \emph{Nilpotent {L}ie groups: structure and applications to
  analysis}, Springer-Verlag, Berlin, 1976, Lecture Notes in Mathematics, Vol.
  562. \MR{56 \#537}

\bibitem{Gross75}
Leonard Gross, \emph{Logarithmic {S}obolev inequalities}, Amer. J. Math.
  \textbf{97} (1975), no.~4, 1061--1083. \MR{54 \#8263}

\bibitem{Gross92}
\bysame, \emph{Logarithmic {S}obolev inequalities and contractivity properties
  of semigroups}, Dirichlet forms (Varenna, 1992), Lecture Notes in Math., vol.
  1563, Springer, Berlin, 1993, pp.~54--88. \MR{95h:47061}

\bibitem{Hormander67}
Lars H{\"o}rmander, \emph{Hypoelliptic second order differential equations},
  Acta Math. \textbf{119} (1967), 147--171. \MR{36 \#5526}

\bibitem{Hsu03}
Elton~P. Hsu, \emph{Stochastic analysis on manifolds}, Graduate Studies in
  Mathematics, vol.~38, American Mathematical Society, Providence, RI, 2002.
  \MR{2003c:58026}

\bibitem{IkedaWatanabe89}
Nobuyuki Ikeda and Shinzo Watanabe, \emph{Stochastic differential equations and
  diffusion processes}, second ed., North-Holland Mathematical Library,
  vol.~24, North-Holland Publishing Co., Amsterdam, 1989. \MR{90m:60069}

\bibitem{KS2}
Shigeo Kusuoka and Daniel Stroock, \emph{Applications of the {M}alliavin
  calculus. {II}}, J. Fac. Sci. Univ. Tokyo Sect. IA Math. \textbf{32} (1985),
  no.~1, 1--76. \MR{86k:60100b}

\bibitem{KS3}
\bysame, \emph{Applications of the {M}alliavin calculus. {III}}, J. Fac. Sci.
  Univ. Tokyo Sect. IA Math. \textbf{34} (1987), no.~2, 391--442.
  \MR{89c:60093}

\bibitem{Li06}
Hong-Quan Li, \emph{Estimation optimale du gradient du semi-groupe de la
  chaleur sur le groupe de {H}eisenberg}, J. Funct. Anal. \textbf{236} (2006),
  no.~2, 369--394. \MR{MR2240167}

\bibitem{Malliavin78b}
Paul Malliavin, \emph{{$C\sp{k}$}-hypoellipticity with degeneracy}, Stochastic
  analysis (Proc. Internat. Conf., Northwestern Univ., Evanston, Ill., 1978),
  Academic Press, New York, 1978, pp.~199--214. \MR{80i:58045a}

\bibitem{Malliavin78c}
\bysame, \emph{{$C\sp{k}$}-hypoellipticity with degeneracy. {II}}, Stochastic
  analysis (Proc. Internat. Conf., Northwestern Univ., Evanston, Ill., 1978),
  Academic Press, New York, 1978, pp.~327--340. \MR{80i:58045b}

\bibitem{Melcher3}
Tai Melcher, \emph{Malliavin calculus for {L}ie group-valued {W}iener
  functions}, Preprint, http://faculty.virginia.edu/melcher (2007).

\bibitem{app}
\bysame, \emph{Some convergence arguments for matrix group-valued {SDE}
  solutions}, http://faculty.virginia.edu/melcher, 2007.

\bibitem{Nualart95}
David Nualart, \emph{The {M}alliavin calculus and related topics}, Probability
  and its Applications (New York), Springer-Verlag, New York, 1995.
  \MR{96k:60130}

\bibitem{Picard}
Jean Picard, \emph{Gradient estimates for some diffusion semigroups}, Probab.
  Theory Related Fields \textbf{122} (2002), no.~4, 593--612. \MR{2003d:58056}

\bibitem{Stein76}
Linda~Preiss Rothschild and E.~M. Stein, \emph{Hypoelliptic differential
  operators and nilpotent groups}, Acta Math. \textbf{137} (1976), no.~3-4,
  247--320. \MR{55 \#9171}

\bibitem{Shigekawa80}
Ichiro Shigekawa, \emph{Derivatives of {W}iener functionals and absolute
  continuity of induced measures}, J. Math. Kyoto Univ. \textbf{20} (1980),
  no.~2, 263--289. \MR{83g:60051}

\bibitem{Stroock80a}
Daniel~W. Stroock, \emph{The {M}alliavin calculus and its application to second
  order parabolic differential equations. {I}}, Math. Systems Theory
  \textbf{14} (1981), no.~1, 25--65. \MR{84d:60092a}

\bibitem{Stroock80b}
\bysame, \emph{The {M}alliavin calculus and its application to second order
  parabolic differential equations. {II}}, Math. Systems Theory \textbf{14}
  (1981), no.~2, 141--171. \MR{84d:60092b}

\bibitem{Varadarajan}
V.~S. Varadarajan, \emph{Lie groups, {L}ie algebras, and their
  representations}, Graduate Texts in Mathematics, vol. 102, Springer-Verlag,
  New York, 1984, Reprint of the 1974 edition. \MR{85e:22001}

\bibitem{Varopoulos86}
N.~Th. Varopoulos, \emph{Analysis on nilpotent groups}, J. Funct. Anal.
  \textbf{66} (1986), no.~3, 406--431. \MR{88h:22014}

\bibitem{Varopoulos92}
N.~Th. Varopoulos, L.~Saloff-Coste, and T.~Coulhon, \emph{Analysis and geometry
  on groups}, Cambridge Tracts in Mathematics, vol. 100, Cambridge University
  Press, Cambridge, 1992. \MR{95f:43008}

\end{thebibliography}

\end{document}